\documentclass[11pt]{article}

\usepackage{amsmath}
\usepackage{amssymb}

 \newcommand{\beq}{\begin{equation}}
\newcommand{\eeq}{\end{equation}}

\newcommand{\R}{{ \mathbb R}}

\newcommand{\C}{{\mathbb  C}}
\newcommand{\Cn}{{\mathbb  C\sp n}}

\newcommand{\D}{{\mathbb D}}

\newcommand{\Z}{{\mathbb  Z}}

\newcommand{\E}{{\mathcal E}}

\newcommand{\cO}{{\mathcal O}}

\newcommand{\cJ}{{\cal J}}
\newcommand{\cI}{{\cal I}}

\newcommand{\m}{{\mathfrak m}}
\newcommand{\n}{{\mathfrak n}}

\newcommand{\codim}{{\operatorname{codim}}}

\newtheorem{theo}{Theorem}
\newtheorem{cor}{Corollary}
\newtheorem{lem}{Lemma}

\begin{document}


\vskip1cm
\begin{center}
{\Large\bf A log canonical threshold test}
\end{center}

\begin{center}
{\large Alexander Rashkovskii}
\end{center}

\vskip1cm

\begin{abstract}\noindent
In terms of log canonical threshold, we characterize plurisubharmonic functions with logarithmic asymptotical behaviour.
\medskip
%
%
\end{abstract}

\section{Introduction and statement of results}

Let $u$ be a plurisubharmonic function on a neighborhood of the origin of $\Cn$. Its {\it log canonical threshold} at $0$,
$$ c_u=
\sup\{c>0:e^{-c\,u}\in L^2_{loc}(0)\},$$
is an important characteristic of asymptotical behavior of $u$ at $0$.
The log canonical threshold $c(\cI)$ of a local ideal in $\cI\subset\cO_0$  can be defined as $c_u$ for the function $u=\log|F|$, where $F=(F_1,\ldots,F_p)$ with $\{F_j\}$ generators of $\cI$. (Surprisingly, the latter notion was introduced later than its plurisubharmonic counterpart.) For general results on log canonical thresholds, including their  computation and applications, we refer to \cite{DK}, \cite{Laz}, \cite{Mu}.

A classical result due to Skoda \cite{Sk} states that
\begin{equation}\label{skoda} c_u\ge \nu_u^{-1},
\end{equation}
where $\nu_u$ is the Lelong number of $u$ at $0$. A more recent result is due to Demailly \cite{D8}: if $0$ is an isolated point of $u^{-1}(-\infty)$, then
\begin{equation}\label{dem} c_u\ge F_n(u):=n\,{e_n(u)^{-1/n}}.\end{equation}
Here $e_k(u)=(dd^cu)^k\wedge (dd^c\log|z|)^{n-k}(0)$ are the Lelong numbers of the currents $(dd^cu)^k$ at $0$ for $k=1,\ldots,n$, and $d=\partial + \bar\partial$, $d^c= (\partial -\bar\partial)/2\pi i$; note that $e_1(u)=\nu_u$. This was extended by Zeriahi \cite{Ze} to all plurisubharmonic functions with well-defined Monge-Amp\`ere operator near $0$.

In \cite{R10}, inequality (\ref{dem}) was used to obtain the `intermediate' bounds
\beq\label{eq:rel} c_u\ge F_k(u):=k\,e_k(u)^{-1/k},\quad 1\le k\le l,\eeq
$l$ being the codimension of an analytic set $A$ containing the unbounded locus $L(u)$ of $u$.
None of the bounds for different values of $k$ can be deduced from the others.

It is worth mentioning that relation (\ref{dem}) was proved in \cite{D8} on the base of a corresponding result for ideals\footnote{A direct proof was given later in \cite{ACKPZ}.} obtained in \cite{dFEM}:
\begin{equation}\label{dFEM} c({\cI})\ge {n}\,{e(\cI)^{-1/n}},\end{equation}
where $e(\cI)$ is the Hilbert-Samuel multiplicity of the (zero-dimensional) ideal $\cI$.  Furthermore, it was shown in \cite{dFEM} that an equality in (\ref{dFEM}) holds if and only if the integral closure of $\cI$ is a power of the maximal ideal $\m_0$. Accordingly, the question of equality in (\ref{dem}) has been raised in \cite{D8} where it was conjectured that, similarly to the case of ideals, the extremal functions would be those with logarithmic singularity at $0$.

The conjecture was proved in \cite{R11} where it was shown that
 \beq\label{eqn} c_u= F_n(u)\eeq
 if and only if the {\it greenification} $g_u$ of $u$ has the asymptotics $g_u(z)=e_1(u)\log|z|+O(1)$ as $z\to 0$. Here the function $g_u$ is the upper semicontinuous regularization of the upper envelope of all negative plurisubharmonic functions $v$ on a bounded neighborhood $D$ of $0$, such that $v\le u+O(1)$ near $0$, see \cite{R7}. Note that if $u=\log|F|$, then $g_u=u+O(1)$ \cite[Prop. 5.1]{RSig2}.

The equality situation in (\ref{skoda}) (i.e., in (\ref{eq:rel}) with $k=1$) was first treated in \cite{BM} and \cite{FaJ2} for the dimension $n=2$: the functions satisfying $c_u=\nu_u^{-1}$ were proved in that case to be of the form $u=c\log|f|+v$, where $f$ is an analytic function, regular at $0$, and $v$ is a plurisubharmonic function with zero Lelong number at $0$. In a recent preprint \cite{QGXZ}, the result was extended to any $n$. This was achieved by a careful slicing technique reducing the general case to the aforementioned two-dimensional result. In addition, it used a regularization result for plurisubharmonic functions with keeping the log canonical threshold (see Lemma~\ref{QGXZ} below).

Concerning inequalities (\ref{eq:rel}), it was shown in \cite{R10} that the only multi-circled plurisubharmonic functions $u(z)=u(|z_1|,\ldots,|z_n|)$ satisfying $c_u=F_l(u)$ are essentially of the form $c\max_{j\in J}\log|z_j|$ for an $l$-tuple $J\subset\{1,\ldots,n\}$. Here we address the question on equalities in the bounds (\ref{eq:rel}) in the general case.

\medskip

We present an approach that is different from that of \cite{QGXZ} and which actually works also for the `intermediate' equality situations. It is based on a recent result of Demailly and Pham Hoang Hiep \cite{DH}: if the complex Monge-Amp\`ere operator $(dd^cu)^n$ is well defined near $0$ and $e_1(u)>0$, then
$$ c_u\ge E_n(u):=\sum_{1\le j\le n}\frac{e_{j-1}(u)}{e_{j}(u)},$$
where $e_0(u)=1$. In particular, this implies (\ref{dem}) and sharpens, for the case of functions with well-defined Monge-Amp\`ere operator, inequality (\ref{skoda}).
Moreover, it is this bound that was used in \cite{R11} to prove the conjecture from \cite{D8} on functions satisfying (\ref{eqn}).

\medskip

Given $1<l\le n$, let $\E_l$  be the collection of all plurisubharmonic functions $u$ whose unbounded loci $L(u)$ have zero $2(n-l+1)$-dimensional Hausdorff measure. For such a function $u$, the currents $(dd^c u)^k$ are well defined for all $k\le l$ \cite{FoSi}. In particular, $u\in\E_l$ if $L(u)$ lies in an analytic variety of codimension at least $l$. Furthermore, we set $\E_1$ to be just the collection of all plurisubharmonic functions near $0$.

Let $c_u(z)$ denote the log canonical threshold of $u$ at $z$ and, similarly, let $e_k(u,z)$ denote the Lelong number of $(dd^cu)^k$ at $z$; in our notation, $c_u(0)=c_u$ and $e_k(u,0)=e_k(u)$. As is known, the sets $\{z:\: c_u(z)\le c\}$ are analytic for all $c>0$.
Our first result describes, in particular, regularity of such a set for $c=c_u$, provided $c_u=F_l(u)$.

For $u\in \E_l$ we set
$$ E_k(u)=\sum_{1\le j\le k}\frac{e_{j-1}(u)}{e_{j}(u)},\quad k\le l.$$

\begin{theo}\label{theo1}  Let $u\in \E_l$ for some $l\ge 1$, and let $e_1(u)>0$. Then
\begin{enumerate}
\item[(i)] $c_u\ge E_k(u)$ for all $k\le l$;
\item[(ii)] $c_u\ge F_k(u)$ for all $k\le l$;
\item[(iii)]  if $u$ satisfies $c_u= F_k(u)$ for some $k\le l$, then $k=l$ and there is a neighborhood $V$ of the origin such that the set $A=\{z:\: c_u(z)\le c_u\}$ is an $l$-codimensional manifold in $V$. Furthermore, $A=\{z:\: e_l(u,z)\ge e_l(u)\}$.
    \end{enumerate}
\end{theo}

For $l=1$, assertion (iii) re-proves the aforementioned result from \cite{QGXZ}. Let $A=\{z_1=0\}$, then the function $u-c_u\log|z_1|$ is locally bounded from above near $A$ and thus extends to a plurisubharmonic function $v$; evidently, $\nu_v=0$. On the other hand, all the functions $u= c_u\log|z_1|+v$ with $\nu_v=0$ satisfy $c_u=\nu_u$.

\medskip
When $l>1$, there are functions $u$ such that  $\{z:\: c_u(z)\le c_u\}$ is an $l$-codimensional manifold, but $c_u>F_l(u)$.
Indeed, let us take $u(z_1,z_2,z_3)=\max\{\log|z_1|, 2\log|z_2|\}\in \E_2$. Then $A=\{z\in\Cn:\: c_u(z)\le c_u\}=\{z_1=z_2=0\}$, while $F_2(u)= \sqrt2<3/2=c_u$. (Note that $c_u=E_2(u)$ in this case.)

Furthermore, the same example shows that the equality $(dd^cu)^2=\delta^2\,[z_1=z_2=0]$ does not imply $u=\delta\log|(z_1,z_2)|+v$ with plurisubharmonic $v$ and $\nu_v=0$.

\medskip

Therefore, in the higher dimensional situation we need to deduce a more precise information on asymptotical behavior of $u$ near $A$.
By analogy with the case $l=n$, it is tempting to make the following conjecture.

\medskip\noindent
{\sl Let $u\in\E_l$, then
\beq\label{eql} c_u= F_l(u)\eeq
 if and only if, for a choice of coordinates $z=(z',z'')\in\C^l\times\C^{n-l}$, the greenification $g_u$ of $u$ near $0$ satisfies
$$g_u=e_1(u)\log|z'|+O(1)\ {\rm as\ } z\to0.$$
}

\medskip
The 'if' direction is obvious in view of $c_u=c_{g_u}$ \cite{R11} and the trivial fact $c_{\log|z'|}=l$, however the reverse statement might be difficult to prove even in the case $l=1$ because that would imply non-existence of a plurisubharmonic function $\phi$ with $e_1(\phi)=0$ and $e_n(\phi)>0$, which is a known open problem. Namely, let such a function $\phi$ exist, and set $u=\phi+\log|z_1|$. Then $1=\nu_u\le c_u\le c_{\log|z_1|}=1$. On the other hand,
for $D=\D^n$, $g_u=g_\phi+\log|z_1|$ and the relation $e_n(\phi)>0$ implies $g_\phi\neq 0$ and thus $\liminf (g_u-\log|z_1|)=-\infty$ when $z\to 0$.

\medskip
What we can prove is the following, slightly weaker statement.

\begin{theo}\label{theo2} If $u\in\E_l$ satisfies (\ref{eql}), then  $e_k(u)=e_1(u)^k$ for all $k\le l$ and, for a choice of coordinates $z=(z',z'')\in\C^l\times\C^{n-l}$, the function $u$ satisfies $u\le e_1(u)\log|z'|+O(1)$ near $0$, while the greenification $g_{u_N}$ of $u_N=\max\{u,N\log|z|\}$ with any $N\ge e_1(u)$ satisfies
\beq\label{gru} g_{u_N}=\max\{e_1(u)\log|z'|,N\log|z''|\}+O(1),\quad z\to0.\eeq
\end{theo}

\medskip

Let us fix a neighborhood $D\subset V$ of $0$ to be the product of unit balls in $\C^l$ and $\C^{n-l}$ and consider the greenifications with respect to $D$. Then the functions $g_{u_N}$ are equal to $\max \{e_1(u)\log|z'|, N\log|z''|\}$ and they converge, as $N\to\infty$, to $e_1(u)\log|z'|\ge g_u$.

Denote, for any bounded neighborhood $D$ of $0$ and any $u$ plurisubharmonic in $D$, $$\tilde g_u=\lim_{N\to\infty} g_{u_N}.$$
where $u_N=\max\{u,N\log|z|\}$. Evidently, $\tilde g_{u}\ge g_u$.

\begin{theo}\label{theo3} Let $u\in\E_l$ be such that $\tilde g_u= g_u$. Then it satisfies (\ref{eql}) if and only if, for a choice of coordinates $z=(z',z'')\in\C^l\times\C^{n-l}$, $g_{u}=e_1(u)\log|z'|+O(1)$ as $z\to0$.

In particular, this is true for $u=\alpha\log|F|+O(1)$, where $F$ is a holomorphic mapping, $F(0)=0$. Moreover, in this case we also have $u=e_1(u)\log|z'|+O(1)$.
\end{theo}

The statement on $\alpha\log|F|$ can be reformulated in algebraic terms as follows. Let $\cI$ be an ideal of the local ring $\cO_0$, and let $V(\cI)$ be its variety: $V(\cI)=\{z:\: f(z)=0\ \forall f\in\cI\}$. If  $\codim_0V(\cI)\ge k$, then the mixed Rees' multiplicity $e_k(\cI,\m_0)$ of $k$ copies of $\cI$ and $n-k$ copies of the maximal ideal $\m_0$ is well defined \cite{BA}. If $k=n$, then, as shown in \cite{D8}, the Hilbert-Samuel multiplicity $e(\cI)$ of $\cI$ equals $e_n(u)$, where, as before, $u=\log|F|$ for generators $\{F_p\}$ of $\cI$. By the polarization formula, $e_k(\cI,\m_0)=e_k(u)$ for all $k$; by a limit transition, this holds true for all $k\le l$ if $\codim_0V(\cI)=l$.

 Bounds (\ref{eq:rel}) specify for this case as
$$ c(\cI)\ge k\,e_k(\cI,\m_0)^{-1/k},\quad 1\le k\le l;$$
from Theorems \ref{theo1} and \ref{theo3} we thus derive

\begin{cor} If $\codim_0V(\cI)= l$ and $ c({\cI})={k}\,{e_k(\cI,\m_0)^{-1/k}}$ for some $k\le l$,
then $k=l$, $V(\cI)$ is an $l$-codimensional hypersurface, regular at $0$, and there exists an ideal $\n_0$ generated by coordinate (smooth transversal) germs $f_1,\ldots,f_l\in\cO_0$ such that $\overline\cI=\n_0^s$ for some $s\in\Z_+$.
\end{cor}

\section{Proofs}

In what follows, we will use the mentioned regularization result by Qi'an Guan and Xiangyu Zhou. Note that its proof rests on the strong openness conjecture from \cite{DK}, proved in \cite{QGXZ1} and \cite{QGXZ2}, see also \cite{Be7}.

\begin{lem} \label{QGXZ} {\rm \cite[Prop. 2.1]{QGXZ}} Let $u$ be a plurisubharmonic function near the origin, $\sigma_u=1$. Then there exists a plurisubharmonic function $\tilde u\ge u$ on a neighborhood of $0$ such that $e^{-2u}-e^{-2\tilde u}$ is integrable on $V$ and $\tilde u$ is locally bounded on $V\setminus \{z:\: c_u(z)\le 1\}$.
\end{lem}

We will also refer to the following uniqueness theorem.

\begin{lem} \label{uniq}{\rm (\cite[Lem. 6.3]{R7}\footnote{For the general case of non-isolated singularities, see \cite[Thm. 3.7]{ACCP}} and \cite[Lem. 1.1]{R11})} If $u$ and $v$ are two plurisubharmonic functions with isolated singularity at $0$, such that $u\le v+O(1)$ near $0$ and $e_n(u)=e_n(v)$, then their greenifications coincide.
\end{lem}

\noindent{\it Proof of Theorem 1.}
 Since all the functionals $u\mapsto c_u,\ E_k(u),\ F_k(u)$ are positive homogeneous of degree $-1$, we can always assume $c_u=1$.

Let $\tilde u$ be the function from Lemma~\ref{QGXZ}. Its unbounded locus $L(\tilde u)$ is contained in the analytic variety $A=\{z:\: c_u(z)\le 1\}$. Since $A\subset L(u)$ and $u\in \E_l$, $\codim\, A\ge l$.

For $\tilde u$, statement (i) is proved in \cite[Thm. 1.4]{R11}. Note that the relation $u\le \tilde u$ implies $e_k(u)\ge e_k(\tilde u)$ for all $k\le l$ and thus $E_l(u)\le E_l(\tilde u)$ \cite{DH}. Since $c_u=c_{\tilde u}$, this gives us (i).

Assertion (ii) follows from (i) by the the arithmetic-geometric mean theorem.

To prove (iii), we first note that (i) implies $c_u\ge E_l(u)>E_k(u)\ge F_k(u)$ for any $k<l$, so we cannot have $c_u=F_k(u)$ unless $k=l$.

Next, if the analytic variety $A$ has codimension $m>l$, then $\tilde u\in \E_m$, so $c_u=c_{\tilde u}\ge E_m({\tilde u})>E_l({\tilde u})\ge E_l(u) \ge F_l(u)$, which contradicts the assumption, so $\codim\, A=l$.

Now we prove that $0$ is a regular point of the variety $A$. By Siu's representation formula,
$$(dd^c u)^l=\sum p_j[A_j] +R$$
on a neighborhood $V$ of $0$,
where $p_j>0$, $[A_j]$ are integration currents along $l$-codimensional analytic varieties containing $0$, and $R$ is a closed positive current such that for any $a>0$ the analytic variety $\{z\in V: \nu(R,z)\ge a\}$ has codimension strictly greater than $l$. If $\nu(R,0)>0$, then for almost all points $z\in A$ we have $e_l(u,z) < e_l(u)$. This implies, by (ii), $c_u(z) > c_u$ for all such points $z$, which is impossible.
The same argument shows that the collection $\{A_j\}$ consists of at most one variety and $0$ is its regular point.
{\hfill$\square$\rm \vskip 2mm}

\bigskip\noindent{\it
Proof of Theorem 2.}
By the arithmetic-geometric mean theorem,  the condition $c_u=F_l(u)$ implies, in view of the inequality $c_u\ge E_l(u)$, the relations
$$\frac{e_{k-1}(u)}{e_{k}(u)}=\frac{e_{j-1}(u)}{e_{j}(u)}$$
for any $k,j\le l$, which gives us $e_k(u)=[e_1(u)]^k$ for all $k\le l$.

Since relation (\ref{gru}) for $e_1(u)=0$ is obvious (in this case $g_{u_N}\equiv 0$), we can assume $e_1(u)=1$.

Note that for any $z$, we have $e_k(u,z)\ge [e_1(u,z)]^k$. As follows from the proof of (iii), the relation $c_u=F_l(u)$ implies then, on a neighborhood $V$ of $0$,
$$A\cap V=\{z\in V:\: c_u(z)\le 1\}=\{z\in V:\: F_l(u,z)\le 1\}=\{z\in V:\: e_k(u,z)\ge 1\}$$
for all $k\le l$.
Moreover, we have $e_k(u,z)=e_1(u,z)^k=1$ for almost all $z\in A\cap V$.

Let us choose, according to Theorem~\ref{theo1}, a coordinate system such that $A\cap V=\{z\in V:\: z_k=0,\ 1\le k\le l\}$. Denote $v(z)=\log|z'|$, $z=(z',z'')\in\C^l\times\C^{n-l}$, then
$A\cap V=\{z:\: e_k(u,z)\ge e_k(v,z)\}$, with equalities almost everywhere.

In particular, we have $u(z)\le \log|z-(0,\zeta'')|+C(\zeta'')$ as $z\to (0,\zeta'')$ for all $z\in \Cn$ and $\zeta''\in \C^{n-l}$ that are close enough to $0$. Assuming $u(z)\le 0$ for all $z$ with $\max|z_k|<2$, we get $u(z)\le \log|z-(0,\zeta'')|$ for all $z\in V$ and $\zeta''\in \C^{n-l}$ with $(0,\zeta'')\in V$. By choosing $\zeta''=z''$ this gives us $u(z)\le v(z)$ on $V$.

Let $u_N=\max\{u, N\log|z|\}$ and $v_N=\max\{v, N\log|z|\}$. Then  $u_N\le v_N$, while for $N\ge 1$ we get, by Demailly's comparison theorem for the Lelong numbers \cite{D1},
$$e_n(u_N)\le (dd^c u)^{l}\wedge (dd^c N\log|z|)^{n-l}(0)=N^{n-l}e_l(u) =N^{n-l}=e_n(v_N).$$
 By Lemma~\ref{uniq}, $g_{u_N}=g_{v_N}$.
{\hfill$\square$\rm \vskip 2mm}

\bigskip\noindent{\it
Proof of Theorem 3. } The only part to prove is the one concerning $u=\alpha\log|F|+O(1)$; we assume $\alpha=1$. As follows from Theorem 2, one can choose coordinates such that the zero set $Z_F$ of $F$ is $\{z:\: z'=0\}\cap V\subset\{0\}\times\C^{n-l}$. Observe that for such a function $u$ we have $e_k(u,z) =e_1(u)^k$ for all $z\in Z_F$ near $0$.

Let $\cI$ be the ideal generated by the components of the mapping $F$. Then, as mentioned in Section~1, $e_l(u)$ equals $e_l(\cI,\m_0)$, the mixed multiplicity of $l$ copies of the ideal $\cI$ and $n-l$ copies of the maximal ideal $\m_0$. By \cite[Prop. 2.9]{BA}, $e_l(\cI,\m_0)$ can be computed as the multiplicity $e(\cJ)$ of the ideal $\cJ$ generated by generic functions $\Psi_1,\ldots,\Psi_l\in \cI$ and $\xi_1,\ldots,\xi_{n-l}\in \m_0$. Since $e(\cJ)=e_l(w)$, where $w=\log|\Psi|$, we have
$ e_l(u)=e_l(w)$.

Let now $v=e_1(u)\log|z'|$, $w_N=\max\{w, N\log|z''|\}$, and $v_N=\max\{v, N\log|z''|\}$. Since $w\le \log|F|+O(1)$, we have from Theorem~\ref{theo2} the inequality $w\le v+O(1)$ and thus  $w_N\le v_N+O(1)$. Note that the mapping $\Psi$ satisfies the {\L}ojasiewicz inequality $|\Psi_0(z)|\ge |z'|^M$ near $0$ for some $M>0$. Therefore, for sufficiently big $N$ we have $w_N=w'_N=\max\{w, N\log|z|\}$. Then, as in the proof of Theorem~\ref{theo2}, we compute
$$e_n(w_N)=e_n(w'_N)\le (dd^c w)^{l}\wedge (dd^c N\log|z|)^{n-l}(0)=N^{n-l}e_l(w) =N^{n-l}e_l(u)= e_n(v_N),$$
which, by Lemma~\ref{uniq}, implies $g_{w_N}=g_{v_N}$ for the greenifications on a bounded neighborhood $D$ of $0$.

We can assume $D=\{|z'|<1\}\times\{|z''|<1\}$, then $g_{v_N}=v_N$, while $g_{w_N}\le w_N$ because the latter function is maximal on $D$ and nonnegative on $\partial D$. Letting $N\to\infty$ we get $ w\ge v$.

Since $ w\le u+O(1)$, we have, in particular, $u\ge v+O(1)$, which, in view of Theorem~\ref{theo2}, completes the proof. {\hfill$\square$\rm \vskip 2mm}

\begin{small}

\end{small}

\vskip1cm

Tek/Nat, University of Stavanger, 4036 Stavanger, Norway

\vskip0.1cm

{\sc E-mail}: alexander.rashkovskii@uis.no

\end{document}